\documentclass[a4paper,oneside,reqno,12pt]{report}
\usepackage{amsmath,amssymb}
\usepackage{texdraw}
\usepackage{amscd,amsfonts}
\usepackage{amsthm}
\usepackage[ansinew]{inputenc}
\usepackage{amsbsy}
\usepackage[french]{babel}

\usepackage[dvips]{color}
\usepackage{graphics}
\usepackage{graphicx}

\newtheorem{thm}{Theorem}
\numberwithin{thm}{section}
\newtheorem{prop}[thm]{Proposition}

\theoremstyle{definition}

\newtheorem{ex}[thm]{Example}

\newtheorem{rem}[thm]{Remark}
\newtheorem{rems}[thm]{Remarks}
\newtheorem{nt}[thm]{Notation}

\newtheorem{dfn}{Definition}[section]
\newtheorem{rl}[thm]{Rule}
\renewcommand{\bibname}{References}

    \DeclareMathOperator{\cd}{\cd}

\newcommand{\dl}{\vskip.5cm}
\newcommand{\mrg}{\noindent}
\newcommand{\pr}{\dl\mrg\textit{ Proof. }}

\renewcommand{\bibname}{References}
\begin{document}
\begin{center}
\textbf{\Large TENSOR PERMUTATION MATRICES IN FINITE DIMENSIONS}\\
\end{center}

\noindent RAKOTONIRINA Christian\\
\noindent {\emph{Institut Supérieur de Technologie d'Antananarivo,
IST-T, BP 8122,\\ Madagascar}}\\ E-mail:
pierrekoto@univ-antananarivo.mg\\

\begin{center}
\textbf{Abstract}
\end{center}
We have generalized the properties with the tensor product, of one
$4\times4$ matrix which is a permutation matrix, and we call a
tensor commutation matrix $2\otimes2$. Tensor commutation matrices
$n \otimes p$ can be construct with or without calculus. A formula
allows us to construct a tensor permutation matrix, which is a
generalization of tensor commutation matrix, has been established.
The expression of an element of a tensor commutation matrix $n
\otimes p$, with $n$ = $p$ , has been generalized in the case of any
element of a tensor permutation matrix. The tensor commutation
matrix $3\otimes3$ has been expressed by using the Gell-Mann
\nolinebreak matrices.\newline

\subsection*{INTRODUCTION}

When we had worked on RAOELINA ANDRIAMBOLOLONA idea on the using
tensor product in Dirac equation \cite{rak03}, \cite{wan01} we had
met the unitary
matrix\\
\begin{center}
$[U_{2\otimes2}]$\;=\;$\left[%
\begin{array}{cccc}
  1 & 0 & 0 & 0 \\
  0 & 0 & 1 & 0 \\
  0 & 1 & 0 & 0 \\
  0 & 0 & 0 & 1 \\
\end{array}%
\right]$
\end{center}
which has the following properties: for any unicolumns and two rows
matrices
\begin{center}$[\alpha]$\;=\;$\left[%
\begin{array}{c}
  \alpha^1 \\
  \alpha^2 \\
\end{array}%
\right]$\;$\in$\;$\mathcal{M}_{2\times1}$($\mathbb{K}$),
$[\beta]$\;=\;$\left[%
\begin{array}{c}
  \beta^1 \\
  \beta^2 \\
\end{array}%
\right]$\;$\in$\;$\mathcal{M}_{2\times1}$($\mathbb{K}$)\\
$[U_{2\otimes2}]$$\cdot$(\;$[\alpha]$$\otimes$$[\beta]$\;)\;=\;$[\beta]$$\otimes$$[\alpha]$
\end{center}
and \cite{rak03} for any two $2\times2$ square matrices, $[A]$,
$[B]$ $\in$$\mathcal{M}_{2\times2}$($\mathbb{K}$)
\begin{center}$[U_{2\otimes2}]$$\cdot$(\;$[A]$$\otimes$$[B]$\;) = (\;$[B]$$\otimes$$[A]$\;)$\cdot$$[U_{2\otimes2}]$
\end{center}
This matrix is frequently found in quantum information theory
\cite{fuj01}, \cite{fad95}, \cite{fra02} where one write, by using
the Pauli matrices \cite{fuj01},\cite{fad95},\cite{fra02}
\begin{center}\begin{equation}\label{e1}$$$[U_{2\otimes2}]$\;=\;$\frac{1}{2}$(\;$\sigma^0$$\otimes$$\sigma^0$
+ $\sigma^i$$\otimes$$\sigma^i$\;)$$
\end{equation}\end{center}
where,
\begin{center}$\sigma^1$ = $\left[%
\begin{array}{cc}
  0 & 1 \\
  1 & 0 \\
\end{array}%
\right]$, $\sigma^2$\;=\;$\left[%
\begin{array}{cc}
  0 & -i \\
  i & 0 \\
\end{array}%
\right]$, $\sigma^3$\;=\;$\left[%
\begin{array}{cc}
  1 & 0 \\
  0 & -1 \\
\end{array}%
\right]$
\end{center}
are the Pauli matrices, $\sigma^0$\;=\;$\left[%
\begin{array}{cc}
  1 & 0 \\
  0 & 1 \\
\end{array}%
\right]$ and according to Einstein convention repeated index denotes
summation over the index in question. We call this matrix a tensor
commutation matrix $2\times2$. A tensor commutation matrix
$3\times3$ has been written by KAZUYUKI FUJII \cite{fuj01} by the
following way\\

\begin{center}$[U_{3\otimes3}]$\;=\;$\left[%
\begin{array}{ccc}
  \left[%
\begin{array}{ccc}
  1 & 0 & 0 \\
  0 & 0 & 0 \\
  0 & 0 & 0 \\
\end{array}%
\right] & \left[%
\begin{array}{ccc}
  0 & 0 & 0 \\
  1 & 0 & 0 \\
  0 & 0 & 0 \\
\end{array}%
\right] & \left[%
\begin{array}{ccc}
  0 & 0 & 0 \\
  0 & 0 & 0 \\
  1 & 0 & 0 \\
\end{array}%
\right] \\
  \left[%
\begin{array}{ccc}
  0 & 1 & 0 \\
  0 & 0 & 0 \\
  0 & 0 & 0 \\
\end{array}%
\right] & \left[%
\begin{array}{ccc}
  0 & 0 & 0 \\
  0 & 1 & 0 \\
  0 & 0 & 0 \\
\end{array}%
\right] & \left[%
\begin{array}{ccc}
  0 & 0 & 0 \\
  0 & 0 & 0 \\
  0 & 1 & 0 \\
\end{array}%
\right] \\
  \left[%
\begin{array}{ccc}
  0 & 0 & 1 \\
  0 & 0 & 0 \\
  0 & 0 & 0 \\
\end{array}%
\right] & \left[%
\begin{array}{ccc}
  0 & 0 & 0 \\
  0 & 0 & 1 \\
  0 & 0& 0 \\
\end{array}%
\right] & \left[%
\begin{array}{ccc}
  0 & 0 & 0\\
  0 & 0 & 0 \\
  0 & 0 & 1 \\
\end{array}%
\right] \\
\end{array}%
\right]$
\end{center}
in order to obtain a conjecture of the form of a tensor commutation
matrix $n\otimes n$, for any $n$ $\in$ $\mathbb{N}^\star$.\indent

We have studied the properties of such matrices in the more general
case, according to the RAOELINA ANDRIAMBOLOLONA approach in linear
algebra and multilinear algebra \cite{rao86}. In establishing
firstly, the theorems on linear operators in intrinsic way, that is
independently of the basis, and after that we demonstrate the
analogous theorems for the matrices.\indent

 Define $[U_{n\otimes
p}]$ as the tensor commutation matrix $n\otimes p$, $n$,
$p$\;$\in$\;$\mathbb{N}^\star$, whose elements are 0 or 1.\newline
In this article we have given two manners to construct $[U_{n\otimes
p}]$ for any $n$ and $p$ $\in$ $\mathbb{N}^\star$, and we have
constructed a formula which allows us to construct the tensor
permutation matrix $[U_{n_1\otimes n_2\otimes\ldots\otimes
n_k}(\sigma)]$, (\;$n_1$, $n_2$, \ldots,
$n_k$\;)\;$\in$\;$\mathbb{N}^\star$ and a formula which gives us the
expression of their elements. From \eqref{e1} it is normal to think
to what about the expression of $[U_{3\otimes3}]$ by using Gell-Mann
matrices. But at first we are obliged to talk a bit about the
definitions of the types of matrices, and after that we are going to
expose the properties of tensor product.\indent

Define $[I_n]$ as the $n\times n$ unit matrix. For the vectors and
covectors we have used the RAOELINA ANDRIAMBOLOLONA's
notations\cite{raoram74}, with overlining for the vectors,
$\overline{x}$, and underlining for the covectors,
$\underline{\varphi}$.\indent

Throughout this article $\mathbb{K}$\;=\;$\mathbb{R}$ or
$\mathbb{C}$.
 \section{Matrices. A generalization}
 If the elements of a matrix are considered as the components of a
 second order tensor, we adopt the habitual notation for a matrix,
 without bracket inside. Whereas if the elements of the matrix are,
 for instance, considered as the components of sixth order tensor,
 three times covariant and three times contravariant, then we
 represent the matrix of the following way, for example\\
 \begin{center}
 $[M]$\;=\;$\left[%
\begin{array}{cc}
  \left[%
\begin{array}{cc}
  \left[%
\begin{array}{cc}
  1 & 0 \\
  1 & 1 \\
\end{array}%
\right] & \left[%
\begin{array}{cc}
  1 & 1 \\
  3 & 2 \\
\end{array}%
\right] \\
  \left[%
\begin{array}{cc}
  0 & 0 \\
  0 & 0 \\
\end{array}%
\right] & \left[%
\begin{array}{cc}
  1 & 1 \\
  1 & 1 \\
\end{array}%
\right] \\
\end{array}%
\right] & \left[%
\begin{array}{cc}
  \left[%
\begin{array}{cc}
  1 & 0 \\
  1 & 2 \\
\end{array}%
\right] & \left[%
\begin{array}{cc}
  7 & 8 \\
  9 & 0 \\
\end{array}%
\right] \\
  \left[%
\begin{array}{cc}
  3 & 4 \\
  5 & 6 \\
\end{array}%
\right] & \left[%
\begin{array}{cc}
  9 & 8 \\
  7 & 6 \\
\end{array}%
\right] \\
\end{array}%
\right] \\
  \left[%
\begin{array}{cc}
  \left[%
\begin{array}{cc}
  1 & 1 \\
  1 & 1 \\
\end{array}%
\right] & \left[%
\begin{array}{cc}
  0 & 0 \\
  3 & 2 \\
\end{array}%
\right] \\
  \left[%
\begin{array}{cc}
  4 & 5 \\
  1 & 6 \\
\end{array}%
\right] & \left[%
\begin{array}{cc}
  1 & 7 \\
  8 & 9 \\
\end{array}%
\right] \\
\end{array}%
\right] & \left[%
\begin{array}{cc}
  \left[%
\begin{array}{cc}
  5 & 4 \\
  3 & 2 \\
\end{array}%
\right] & \left[%
\begin{array}{cc}
  1 & 0 \\
  1 & 2 \\
\end{array}%
\right] \\
  \left[%
\begin{array}{cc}
  3 & 4 \\
  5 & 6 \\
\end{array}%
\right] & \left[%
\begin{array}{cc}
  7 & 8 \\
  9 & 0 \\
\end{array}%
\right] \\
\end{array}%
\right] \\
  \left[%
\begin{array}{cc}
  \left[%
\begin{array}{cc}
  1 & 2 \\
  3 & 4 \\
\end{array}%
\right] & \left[%
\begin{array}{cc}
  9 & 8 \\
  7 & 6 \\
\end{array}%
\right] \\
  \left[%
\begin{array}{cc}
  5 & 6 \\
  7 & 8 \\
\end{array}%
\right] & \left[%
\begin{array}{cc}
  5 & 4 \\
  3 & 2 \\
\end{array}%
\right] \\
\end{array}%
\right] & \left[%
\begin{array}{cc}
  \left[%
\begin{array}{cc}
  9 & 8 \\
  7 & 6 \\
\end{array}%
\right] & \left[%
\begin{array}{cc}
  5 & 4 \\
  3 & 2 \\
\end{array}%
\right] \\
  \left[%
\begin{array}{cc}
  1 & 0 \\
  1 & 2 \\
\end{array}%
\right] & \left[%
\begin{array}{cc}
  3 & 4 \\
  5 & 6 \\
\end{array}%
\right] \\
\end{array}%
\right] \\
\end{array}%
\right]$\linebreak
 \end{center}

$[M]$\;=\;$\left(%
  M_{j_1j_2j_3}^{i_1i_2i_3} \\
\right)$\\

$i_1i_2i_3$\;=\;$111,\; 112,\; 121,\; 122,\; 211,\; 212,\; 221,\;
222,\; 311,\;312,\;
321,\; 322$\\
row indices\\

 $j_1j_2j_3$\;=\;$111,\; 112,\; 121,\; 122,\; 211,\; 212,\; 221,\; 222$\\
 column indices\\
 The first indices $i_1$  and $j_1$  are the indices of the outside bracket
 which we call the first order bracket; the second indices $i_2$  and
 $j_2$
 are the indices of the next bracket  which we call the second
 order brackets; the third indices $i_3$  and $j_3$  are the indices of the
 most interior bracket, of this example, which we call third
 order brackets. So, for instance, $M_{121}^{321}$\;=\;5. \\
 If we delete the third order brackets, then the elements of the matrix $[M]$
  are considered as the components of a forth order tensor, twice
  contravariant and twice covariant.\\
  \indent
  Consider a more general case $[M]$\;=\;$\left(%
  M_{j_1j_2\ldots j_k}^{i_1i_2\ldots i_k} \\
\right)$ where the elements of $[M]$ are
  considered as the components of a $2k$-th order tensor, $k$ times contravariant
  and $k$ times covariant. The first order bracket is the bracket of a  $n_1\times\;m_1$-
  dimensional matrix; the second order brackets are the brackets of a $n_2\times\;m_2$-
  dimensional matrices; … ; the $k$-th brackets are the the brackets of the
  $n_k\times\;m_k$-dimensional matrices. $[M]$\;=\;$\left(%
  \gamma_t^s \\
\right)_{1\;\leq\;s\;\leq\;n_1n_2\ldots
n_k,\;1\;\leq\;t\;\leq\;m_1m_2\ldots m_k}$ if the elements of $[M]$
are considered as the components a second order tensor, once
contravariant and once covariant. Then, \cite{rao77}
\numberwithin{equation}{section}
\begin{equation}\label{e12}s=n_kn_{k-1}\ldots n_3n_2(i_1-1)+n_kn_{k-1}\ldots n_3(i_2-1)+\ldots
+n_k(i_{k-1}-1)+i_k
\end{equation}
\begin{equation}\label{e13}t=m_km_{k-1}\ldots m_3m_2(j_1-1)+m_km_{k-1}\ldots
m_3(j_2-1)+\ldots +m_k(j_{k-1}+j_k
\end{equation}
 The elements of the matrix $[N]$\;=\;$(N_k^{ij})$\;=\;$\left[%
\begin{array}{c}
  \left[%
\begin{array}{cc}
  1 & 1 \\
  1 & 1 \\
\end{array}%
\right] \\
  \left[%
\begin{array}{cc}
  0 & 1 \\
  2 & 3 \\
\end{array}%
\right] \\
  \left[%
\begin{array}{cc}
  4 & 5 \\
  6 & 7 \\
\end{array}%
\right] \\
  \left[%
\begin{array}{cc}
  8 & 9 \\
  0 & 1 \\
\end{array}%
\right] \\
\end{array}%
\right]$   , with interior brackets, can
 be considered as the components of a third order tensor, twice contravariant
 and once covariant. Then, for example, $(N_2^{12})$\;=\;1.
 \section{Tensor Product of Matrices}

\begin{dfn}
Consider $[A]$\;=\;$(A^i_j)$\;$\in$\;$\mathcal{M}_{m\times
n}$($\mathbb{K}$), $[B]$\;=\;$(B^i_j)$\;$\in$\;$\mathcal{M}_{p\times
r}$($\mathbb{K}$). The matrix defined by
\begin{center}
$[A]$$\otimes$ $[B]$\;=\;$\left[%
\begin{array}{ccccc}
  $A$^1_1$[B]$ & \ldots & $A$^1_j$[B]$ & \ldots & $A$^1_n$[B]$ \\
  \vdots &  & \vdots &  & \vdots \\
  $A$^i_1$[B]$ & \ldots & $A$^i_j$[B]$ & \ldots & $A$^i_n$[B]$ \\
  \vdots &  & \vdots &  & \vdots \\
  $A$^m_1$[B]$ & \ldots & $A$^m_j$[B]$ & \ldots & $A$^m_n$[B]$ \\
\end{array}%
\right]$
\end{center}
obtained by suppressing the brackets after the multiplications by
scalar, $A^i_j$$[B]$, is called the tensor product of the matrix
$[A]$ by the matrix $[B]$.
\begin{center}
$[A]$$\otimes$ $[B]$\;$\in$\;$\mathcal{M}_{mp\times
nr}$($\mathbb{K}$)
\end{center}
\end{dfn}
\numberwithin{thm}{section}
\begin{thm}
Consider the linear operators $A$\;$\in$\;$\mathcal{L}( \mathcal{E},
\mathcal{F} )$ , $B$\;$\in$\;$\mathcal{L}( \mathcal{G}, \mathcal{H}
)$. $[A]$ is the matrix of A in the couple of basis
$\left(%
  (\overline{e_i})_{1\;\leq i\;\leq n}, (\overline{f_j})_{1\;\leq
j\;\leq m} \\
\right)$, $[B]$  the one of B in
$\left(%
(\overline{g_k})_{1\;\leq
k\;\leq r}, (\overline{h_l})_{1\;\leq l\;\leq p}\\
\right)$. Then, $[A]$$\otimes$ $[B]$ is the matrix of
$A$$\otimes$$B$ in the couple of basis ($\mathcal{B}$,
$\mathcal{B}_1$), where\\
$\mathcal{B}$\;=\;$(\overline{e_1} \otimes \overline{g_1},
\overline{e_1} \otimes \overline{g_2}, \ldots ,\overline{e_1}
\otimes \overline{g_r}, \overline{e_2} \otimes \overline{g_1},
\overline{e_2} \otimes \overline{g_2}, \ldots , \overline{e_2}
\otimes \overline{g_r}, \ldots , \overline{e_n} \otimes
\overline{g_1}, \overline{e_n} \otimes
\overline{g_2}, \ldots , \overline{e_n} \otimes \overline{g_r})$\\
\end{thm}
Define\\
 $\mathcal{B}\;=\;\left(%
 \overline{e_i} \otimes \overline{g_k}\\
\right)_{1\;\leq i\;\leq n, \\
        1\;\leq k\;\leq r}$\;=\;$\left(%
  (\overline{e_i})_{1\;\leq i\;\leq n} \\
\right)$ $\otimes$ $\left(%
(\overline{g_k})_{1\;\leq
k\;\leq r}\\
\right)$\\
 $\mathcal{B}_1\;=\;\left(%
 \overline{f_j}\otimes \overline{h_l}\\
\right)_{1\;\leq j\;\leq m, \\
        1\;\leq l\;\leq p}$\;=\;$\left(%
  (\overline{f_j})_{1\;\leq j\;\leq m} \\
\right)$ $\otimes$ $\left(%
(\overline{h_l})_{1\;\leq
l\;\leq p}\\
\right)$\\

Let
\begin{center}
$[A]$ $\otimes$ $[B]$\;=\;$\left(%
  \gamma_t^s
\right)_{1\;\leq\;s\;\leq\;mp,\;1\;\leq\;t\;\leq\;nr}$
\end{center}
where,
\begin{center}
$\gamma_t^s$\;=\;$A^{i_1}_{j_1}$$B^{i_2}_{j_2}$,
\end{center}
according to the formulas \eqref{e12} and \eqref{e13}
\begin{center}
$s$\;=\;$p(i_1-1)+i_2$\\
$t$\;=\;$r(j_1-1)+j_2$
\end{center}

One can also define elements of tensor product of the matrix $[A]$
by the matrix $[B]$ by the the following way
\begin{center}
$[A]$ $\otimes$
$[B]$\;=\;$(C_{j_1j_2}^{i_1i_2})$\;=\;($A^{i_1}_{j_1}$$B^{i_2}_{j_2}$)
\end{center}
(cf. for example \cite{fuj01} )
 where,\\
  $i_1i_2$ are row indices\\
  $j_1j_2$ are column indices.\\
In the first notation the elements of $[A]$$\otimes$ $[B]$ are
considered as components of a second order tensor , whereas in the
second one, the elements of $[A]$$\otimes$ $[B]$ are considered as
components of a forth order tensor. As an example,
\begin{center}
$[A]$\;=\;$\left[%
\begin{array}{cc}
  1 & 0 \\
  -1 & 2 \\
\end{array}%
\right], \left[%
\begin{array}{ccc}
  0 & 1 & -1 \\
  2 & 3 & 4 \\
\end{array}%
\right]$
\end{center}

\begin{center}
$[A]$ $\otimes$ $[B]$\;=\;$\left[%
\begin{array}{cc}
  \left[%
\begin{array}{ccc}
  0 & 1 & -1 \\
  2 & 3 & 4 \\
\end{array}%
\right] & \left[%
\begin{array}{ccc}
  0 & 0 & 0 \\
  0 & 0 & 0 \\
\end{array}%
\right] \\
  \left[%
\begin{array}{ccc}
  0 & -1 & 1 \\
  -2 & -3 & -4 \\
\end{array}%
\right] & \left[%
\begin{array}{ccc}
  0 & 2 & -2 \\
  4 & 6 & 8 \\
\end{array}%
\right] \\
\end{array}%
\right]$
\end{center}
the result of the calculus for the second notation, whereas
\begin{center}
$[A]$ $\otimes$ $[B]$\;=\;$\left[%
\begin{array}{cccccc}
  0 & 1 & -1 & 0 & 0 & 0 \\
  2 & 3 & 4 & 0 & 0 & 0 \\
  0 & -1 & 1 & 0 & 2 & -2 \\
  -2 & -3 & -4 & 4 & 6 & 8 \\
\end{array}%
\right]$
\end{center}
the result of the calculus for the first notation.
\begin{prop}
$\left(%
[B_1]\cdot[A_1]%
\right)$ $\otimes$ $\left(%
[B_2]\cdot[A_2]
\right)$\;=\;$\left(%
[B_1] \otimes [B_2]
\right)\cdot\left(%
[A_1] \otimes [A_2]
\right)$ for any matrices $[B_1]$, $[A_1]$, $[B_2]$, $[A_2]$ if the
habitual matricial products $[B_1]\cdot[A_1]$ and $[B_2]\cdot[A_2]$
are defined.
\end{prop}
\begin{prop}
$[I_n] \otimes [I_m]$\;=\;$[I_{nm}]$
\end{prop}
\begin{prop}
Tensor product of matrices is associative.
\end{prop}
\begin{thm}\label{thm25}
Consider $\left(%
 [A_i]%
\right)_{1\leq i\leq n\times m}$ a basis of
$\mathcal{M}_{n\times m}$($\mathbb{K}$), $\left(%
[B_j] %
\right)_{1\leq j\leq p\times r}$ a basis of
$\mathcal{M}_{p\times r}$($\mathbb{K}$). Then, $\left(%
 [A_i]\otimes [B_j]%
\right)_{1\leq i\leq n\times m, 1\leq j\leq p\times r}$ is a basis
of $\mathcal{M}_{np\times mr}$($\mathbb{K}$).
\end{thm}
\section{Tensor permutation operators}
\begin{dfn}
Consider $k$  $\mathbb{K}$- vector spaces $\mathcal{E}_1$,
$\mathcal{E}_2$, \ldots, $\mathcal{E}_k$, $\sigma$ a permutation on
\{${1, 2,\ldots , k}$\}. Call  $\sigma$- tensor permutation operator
a linear operator $U_{\sigma}$  from $\mathcal{E}_1$$\otimes$
$\mathcal{E}_2$$\otimes$\ldots$\otimes$ $\mathcal{E}_k$ to
$\mathcal{E}_{\sigma(1)}$$\otimes$
$\mathcal{E}_{\sigma(2)}$$\otimes$\ldots$\otimes$
$\mathcal{E}_{\sigma(k)}$, $U_{\sigma}$$\in$
$\mathcal{L}$($\mathcal{E}_1$$\otimes$
$\mathcal{E}_2$$\otimes$\ldots$\otimes$ $\mathcal{E}_k$,
$\mathcal{E}_{\sigma(1)}$$\otimes$
$\mathcal{E}_{\sigma(2)}$$\otimes$\ldots$\otimes$
$\mathcal{E}_{\sigma(k)}$ ), defined by
\begin{center}
$U_{\sigma}$($\overline{x_1}$\;$\otimes$\;$\overline{x_2}$\;
$\otimes$\;\ldots\;$\otimes$\;$\overline{x_k}$)\;=\;
$\overline{x_{\sigma(1)}}$\;$\otimes$\;$\overline{x_{\sigma(2)}}$\;$\otimes$\;\ldots
\;$\otimes$\;$\overline{x_{\sigma(k)}}$
\end{center}
for all $\overline{x_1}$\;$\in$\;$\mathcal{E}_1$,
$\overline{x_2}$\;$\in$\;$\mathcal{E}_2$, \ldots,
$\overline{x_k}$\;$\in$\;$\mathcal{E}_k$.\\ If $n$\;=\;2, then say
that $U_{\sigma}$  is a tensor commutation operator.
\end{dfn}
\begin{rem}
$U_{\sigma}$ is an isomorphism and ${U_{\sigma}}^{-1}$ is a tensor
permutation operator too.
\end{rem}
\begin{thm}\label{tm32}
Consider k  $\mathbb{K}$- vector spaces $\mathcal{E}_1$,
$\mathcal{E}_2$, \ldots, $\mathcal{E}_k$, $\sigma$ a permutation on
\{${1, 2,\ldots , k}$\}, $U_{\sigma}$\;$\in\;$
$\mathcal{L}$($\mathcal{E}_1$\;$\otimes$\;
$\mathcal{E}_2$\;$\otimes$\;\ldots\;$\otimes$ $\mathcal{E}_k$,
$\mathcal{E}_{\sigma(1)}$\;$\otimes$\;
$\mathcal{E}_{\sigma(2)}$\;$\otimes$\;\ldots\;$\otimes$\;
$\mathcal{E}_{\sigma(k)}$) $\sigma$- tensor permutation operator.
Then, for all \; $\phi_1$\;$\in$\; $\mathcal{L}$($\mathcal{E}_1$),
$\phi_2$\;$\in$\; $\mathcal{L}$($\mathcal{E}_2$), \ldots,
$\phi_k$\;$\in$\;$\mathcal{L}$($\mathcal{E}_k$)
\begin{center}
$U_{\sigma}$$\cdot$($\phi_1$\;$\otimes$\;$\phi_2$\;$\otimes$\;\ldots\;$\otimes$\;$\phi_k$)
=
($\phi_{\sigma(1)}$\;$\otimes$\;$\phi_{\sigma(2)}$\;$\otimes$\;\ldots\;$\otimes$\;
$\phi_{\sigma(k)}$)$\cdot$$U_{\sigma}$
\end{center}
\end{thm}

\pr
$\phi_1$\;$\otimes$\;$\phi_2$\;$\otimes$\;\ldots\;$\otimes$\;$\phi_k$\;$\in\;
$$\mathcal{L}(\mathcal{E}_1\;\otimes\;
\mathcal{E}_2\;\otimes\;\ldots\;\otimes
\mathcal{E}_k)$, thus\\
$U_{\sigma}$$\cdot$($\phi_1$\;$\otimes$\;$\phi_2$\;$\otimes$\ldots$\otimes$\;$\phi_k$)$\in$
$\mathcal{L}$($\mathcal{E}_1$\;$\otimes$\;
$\mathcal{E}_2$\;$\otimes$\;\ldots\;$\otimes$ $\mathcal{E}_k$,
$\mathcal{E}_{\sigma(1)}$\;$\otimes$
$\mathcal{E}_{\sigma(2)}$\;$\otimes$\;\ldots\;$\otimes$\;
$\mathcal{E}_{\sigma(k)}$)
$\phi_{\sigma(1)}$\;$\otimes$\;$\phi_{\sigma(2)}$\;$\otimes$\;\ldots\;$\otimes$\;
$\phi_{\sigma(k)}$\;$\in\;$$\mathcal{L}(\mathcal{E}_{\sigma(1)}\;\otimes\;
\mathcal{E}_{\sigma(2)}\;\otimes\;\ldots\;\otimes\;
\mathcal{E}_{\sigma(k)})$, thus\\
($\phi_{\sigma(1)}$\;$\otimes$\;$\phi_{\sigma(2)}$\;$\otimes$\ldots$\otimes$\;
$\phi_{\sigma(k)}$)$\cdot$$U_{\sigma}$$\in$
$\mathcal{L}$($\mathcal{E}_1$\;$\otimes$\;
$\mathcal{E}_2$\;$\otimes$\ldots$\otimes$ $\mathcal{E}_k$,
$\mathcal{E}_{\sigma(1)}$\;$\otimes$
$\mathcal{E}_{\sigma(2)}$\;$\otimes$\ldots$\otimes$\;
$\mathcal{E}_{\sigma(k)}$)\\
If $\overline{x_1}$\;$\in$\;$\mathcal{E}_1$,
$\overline{x_2}$\;$\in$\;$\mathcal{E}_2$, \ldots,
$\overline{x_k}$\;$\in$\;$\mathcal{E}_k$,\\
$U_{\sigma}$$\cdot$($\phi_1$\;$\otimes$\;$\phi_2$\;$\otimes$\ldots$\otimes$\;$\phi_k$)
($\overline{x_1}$\;$\otimes$\;$\overline{x_2}$\;
$\otimes$\;\ldots\;$\otimes$\;$\overline{x_k}$)\;=\;$U_{\sigma}$[$\phi_1$($\overline{x_1}$)\;
$\otimes$\;$\phi_2$($\overline{x_2}$)\;$\otimes$\ldots$\otimes$\;$\phi_k$($\overline{x_k}$)]\;\\
=\;$\phi_{\sigma(1)}$($\overline{x_{\sigma(1)}}$)\;
$\otimes$\;$\phi_{\sigma(2)}$($\overline{x_{\sigma(2)}}$)\;$\otimes$\ldots$\otimes$\;
$\phi_{\sigma(k)}$($\overline{x_{\sigma(k)}}$)\\
( Since $U_{\sigma}$ is a  tensor permutation operator)\\
=\;($\phi_{\sigma(1)}$\;$\otimes$\;$\phi_{\sigma(2)}$\;$\otimes$\;\ldots\;$\otimes$\;
$\phi_{\sigma(k)}$)($\overline{x_{\sigma(1)}}$\;$\otimes$\;$\overline{x_{\sigma(2)}}$\;$\otimes$\;\ldots
\;$\otimes$\;$\overline{x_{\sigma(k)}}$)\\=\;($\phi_{\sigma(1)}$\;$\otimes$\;
$\phi_{\sigma(2)}$\;$\otimes$\;\ldots\;$\otimes$\;
$\phi_{\sigma(k)}$)$\cdot$$U_{\sigma}$($\overline{x_1}$\;$\otimes$\;$\overline{x_2}$\;
$\otimes$\;\ldots\;$\otimes$\;$\overline{x_k}$). \qed\\

\begin{thm}
If  $U_{\sigma}$ is a $\sigma$- tensor permutation operator, then
its transpose ${U_{\sigma}}^t$  is a ${\sigma}^{-1}$-tensor
permutation operator.\\
\end{thm}
\pr Let us consider k  $\mathbb{K}$- vector spaces $\mathcal{E}_1$,
$\mathcal{E}_2$, \ldots, $\mathcal{E}_k$ of finite dimensions,
$U_{\sigma}$\;$\in\;$ $\mathcal{L}$($\mathcal{E}_1$\;$\otimes$\;
$\mathcal{E}_2$\;$\otimes$\;\ldots\;$\otimes$ $\mathcal{E}_k$,
$\mathcal{E}_{\sigma(1)}$\;$\otimes$\;
$\mathcal{E}_{\sigma(2)}$\;$\otimes$\;\ldots\;$\otimes$\;
$\mathcal{E}_{\sigma(k)}$) $\sigma$- tensor permutation operator.
Then ${U_{\sigma}}^t$\;$\in\;$ $\mathcal{L}$(
${\mathcal{E}_{\sigma(1)}}^{\star}$\;$\otimes$\;
${\mathcal{E}_{\sigma(2)}}^{\star}$\;$\otimes$\;\ldots\;$\otimes$\;
${\mathcal{E}_{\sigma(k)}}^{\star}$,
${\mathcal{E}_1}^{\star}$\;$\otimes$\;
${\mathcal{E}_2}^{\star}$\;$\otimes$\;\ldots\;$\otimes$
${\mathcal{E}_k}^{\star}$).\\
Let
$\underline{\varphi^{\sigma(1)}}$\;$\in\;$${\mathcal{E}_{\sigma(1)}}^{\star}$,
 $\underline{\varphi^{\sigma(2)}}$\;$\in\;$${\mathcal{E}_{\sigma(2)}}^{\star}$,
 \ldots$\underline{\varphi^{\sigma(k)}}$\;$\in\;$${\mathcal{E}_{\sigma(k)}}^{\star}$,
 $\overline{x_1}$\;$\in$\;$\mathcal{E}_1$,
$\overline{x_2}$\;$\in$\;$\mathcal{E}_2$, \ldots,
$\overline{x_k}$\;$\in$\;$\mathcal{E}_k$.\\
${U_{\sigma}}^t\cdot\left(
\underline{\varphi^{\sigma(1)}}\;\otimes\;\underline{\varphi^{\sigma(2)}}
\;\otimes\ldots\otimes\;\underline{\varphi^{\sigma(k)}}
\right)$($\overline{x_1}$\;$\otimes$\;$\overline{x_2}$\;
$\otimes$\;\ldots\;$\otimes$\;$\overline{x_k}$)\;\\

=\;$\underline{\varphi^{\sigma(1)}}\;\otimes\;\underline{\varphi^{\sigma(2)}}
\;\otimes\ldots\otimes\;\underline{\varphi^{\sigma(k)}}$\;[$U_{\sigma}$($\overline{x_1}$\;$\otimes$\;$\overline{x_2}$\;
$\otimes$\;\ldots\;$\otimes$\;$\overline{x_k}$)]\\
( by definition of operator transpose \cite{rao86})\\
=\;$\underline{\varphi^{\sigma(1)}}\;\otimes\;\underline{\varphi^{\sigma(2)}}
\;\otimes\ldots\otimes\;\underline{\varphi^{\sigma(k)}}$($\overline{x_{\sigma(1)}}$\;$\otimes$\;$\overline{x_{\sigma(2)}}$\;$\otimes$\;\ldots
\;$\otimes$\;$\overline{x_{\sigma(k)}}$)\\

=\;$\underline{\varphi^{\sigma(1)}}(\overline{x_{\sigma(1)}})\otimes\underline{\varphi^{\sigma(2)}}(\overline{x_{\sigma(2)}})
\otimes\ldots\otimes\underline{\varphi^{\sigma(k)}}(\overline{x_{\sigma(k)}})$
=\;$\underline{\varphi^{\sigma(1)}}(\overline{x_{\sigma(1)}})\underline{\varphi^{\sigma(2)}}(\overline{x_{\sigma(2)}})
\ldots\underline{\varphi^{\sigma(k)}}(\overline{x_{\sigma(k)}})$\\

(because $\underline{\varphi^{\sigma(i)}}(\overline{x_{\sigma(i)}})$
are
elements of $\mathbb{K}$)\\

=\;$\underline\varphi^1(\overline{x_1})\otimes
\underline\varphi^2(\overline{x_2})
\otimes \ldots\otimes \underline\varphi^k(\overline{x_k})$\\

=\;$(\underline\varphi^1\otimes \underline\varphi^2\otimes
\ldots\otimes
\underline\varphi^k)(\overline{x_1}\;\otimes\;\overline{x_2}\;
\otimes\;\ldots\;\otimes\;\overline{x_k})$\\
We have\\
${U_{\sigma}}^t\left(
\underline{\varphi^{\sigma(1)}}\;\otimes\;\underline{\varphi^{\sigma(2)}}
\;\otimes\ldots\otimes\;\underline{\varphi^{\sigma(k)}}
\right)$($\overline{x_1}$\;$\otimes$\;$\overline{x_2}$\;
$\otimes$\;\ldots\;$\otimes$\;$\overline{x_k}$)\;\\

=\;$(\underline\varphi^1\otimes \underline\varphi^2\otimes
\ldots\otimes
\underline\varphi^k)(\overline{x_1}\;\otimes\;\overline{x_2}\;
\otimes\;\ldots\;\otimes\;\overline{x_k})$\\
for all $\overline{x_1}$\;$\in$\;$\mathcal{E}_1$,
$\overline{x_2}$\;$\in$\;$\mathcal{E}_2$, \ldots,
$\overline{x_k}$\;$\in$\;$\mathcal{E}_k$.\\
Hence,\\
${U_{\sigma}}^t\cdot\left(
\underline{\varphi^{\sigma(1)}}\;\otimes\;\underline{\varphi^{\sigma(2)}}
\;\otimes\ldots\otimes\;\underline{\varphi^{\sigma(k)}}
\right)$=\;$(\underline\varphi^1\otimes \underline\varphi^2\otimes
\ldots\otimes \underline\varphi^k)$ \qed
\section{Tensor permutation matrices}
\begin{dfn}Let us consider $k$ finite dimensional $\mathbb{K}$- vector spaces
$\mathcal{E}_1$, $\mathcal{E}_2$, \ldots, $\mathcal{E}_k$ of
dimensions $n_1$, $n_2$, \ldots, $n_k$ respectively,
$U_{\sigma}$\;$\in\;$ $\mathcal{L}$($\mathcal{E}_1$\;$\otimes$\;
$\mathcal{E}_2$\;$\otimes$\;\ldots\;$\otimes$ $\mathcal{E}_k$,
$\mathcal{E}_{\sigma(1)}$\;$\otimes$\;
$\mathcal{E}_{\sigma(2)}$\;$\otimes$\;\ldots\;$\otimes$\;
$\mathcal{E}_{\sigma(k)}$) $\sigma$- tensor permutation operator.
Let\\
$\mathcal{B}_1$\;=\;$(\overline{e_{11}}\;\otimes\;\overline{e_{12}}\;
\otimes\;\ldots\;\otimes\;\overline{e_{1n_1}})$ be a basis of
$\mathcal{E}_1$\\
$\mathcal{B}_2$\;=\;$(\overline{e_{21}}\;\otimes\;\overline{e_{22}}\;
\otimes\;\ldots\;\otimes\;\overline{e_{2n_2}})$ be a basis of
$\mathcal{E}_2$\\\ldots\\
$\mathcal{B}_k$\;=\;$(\overline{e_{k1}}\;\otimes\;\overline{e_{k2}}\;
\otimes\;\ldots\;\otimes\;\overline{e_{kn_k}})$ be a basis of
$\mathcal{E}_k$\\
$[U_{\sigma}]$ the matrix of $U_{\sigma}$ with respect to a couple
of basis ($\mathcal{B}_1$$\otimes$$\mathcal{B}_2$$\otimes$
\ldots$\otimes$ $\mathcal{B}_k$, $\mathcal{B}_{\sigma(1)}$$\otimes$
$\mathcal{B}_{\sigma(2)}$$\otimes$ \ldots$\otimes$
$\mathcal{B}_{\sigma(k)}$ ). $[U_{\sigma}]$ is a square matrix of
dimension $n_1\times n_2\times\ldots \times n_k$ independent of the
basis $\mathcal{B}_1$, $\mathcal{B}_2$,\ldots , $\mathcal{B}_k$.
Call this matrix $\sigma$-tensor permutation matrix $n_1\otimes
n_2\otimes\ldots \otimes n_k$.
\end{dfn}
\begin{thm}\label{tm41}
$[U_{\sigma}]$ a $\sigma$-tensor permutation matrix $n_1\otimes
n_2\otimes\ldots \otimes n_k$ if only if, for all
$[\alpha_1]\;\in\;\mathcal{M}_{n_1\times 1}(\mathbb{K})$,
$[\alpha_2]\;\in\;\mathcal{M}_{n_2\times 1}(\mathbb{K})$,\ldots,
$[\alpha_k]\;\in\;\mathcal{M}_{n_k\times 1}(\mathbb{K})$
\begin{center}
$[U_{\sigma}]\cdot([\alpha_1]\otimes [\alpha_2]\otimes \ldots\otimes
[\alpha_k])$\;=\;$[\alpha_{\sigma(1)}]\otimes
[\alpha_{\sigma(2)}]\otimes \ldots\otimes [\alpha_{\sigma(k)}]$
\end{center}
\end{thm}
\pr $"\Longrightarrow"$ Let
$\overline{\alpha_1}$\;$\in$\;$\mathcal{E}_1$,
$\overline{\alpha_2}$\;$\in$\;$\mathcal{E}_2$, \ldots,
$\overline{\alpha_k}$\;$\in$\;$\mathcal{E}_k$, $[\alpha_1]$,
$[\alpha_2]$,\ldots, $[\alpha_k]$ unicolumn matrices formed by the
components of $\overline{\alpha_1}$, $\overline{\alpha_2}$,\ldots,
$\overline{\alpha_k}$ respectively with respect to $\mathcal{B}_1$,
$\mathcal{B}_2$,\ldots , $\mathcal{B}_k$.
\begin{center}
$U_{\sigma}$($\overline{\alpha_1}$\;$\otimes$\;$\overline{\alpha_2}$\;
$\otimes$\;\ldots\;$\otimes$\;$\overline{\alpha_k}$)\;=\;
$\overline{\alpha_{\sigma(1)}}$\;$\otimes$\;$\overline{\alpha_{\sigma(2)}}$\;$\otimes$\;\ldots
\;$\otimes$\;$\overline{\alpha_{\sigma(k)}}$
\end{center}
The components of
$\overline{\alpha_1}$\;$\otimes$\;$\overline{\alpha_2}$\;
$\otimes$\;\ldots\;$\otimes$\;$\overline{\alpha_k}$ with respect to
$\mathcal{B}_1$$\otimes$$\mathcal{B}_2$$\otimes$ \ldots$\otimes$
$\mathcal{B}_k$ are the elements of the unicolumn matrix
$[\alpha_{1}]\otimes [\alpha_{2}]\otimes \ldots\otimes [\alpha_{k}]$
and the
$\overline{\alpha_{\sigma(1)}}$\;$\otimes$\;$\overline{\alpha_{\sigma(2)}}$\;$\otimes$\;\ldots
\;$\otimes$\;$\overline{\alpha_{\sigma(k)}}$ ones, with respect to
$\mathcal{B}_{\sigma(1)}$$\otimes$
$\mathcal{B}_{\sigma(2)}$$\otimes$ \ldots$\otimes$
$\mathcal{B}_{\sigma(k)}$, are the elements of the unicolumn matrix
$[\alpha_{\sigma(1)}]\otimes [\alpha_{\sigma(2)}]\otimes
\ldots\otimes [\alpha_{\sigma(k)}]$.\\
Hence,
\begin{center}
$[U_{\sigma}]\cdot([\alpha_1]\otimes [\alpha_2]\otimes \ldots\otimes
[\alpha_k])$\;=\;$[\alpha_{\sigma(1)}]\otimes
[\alpha_{\sigma(2)}]\otimes \ldots\otimes [\alpha_{\sigma(k)}]$
\end{center}
$"\Longleftarrow"$ Suppose that for all
$[\alpha_1]\;\in\;\mathcal{M}_{n_1\times 1}(\mathbb{K})$,
$[\alpha_2]\;\in\;\mathcal{M}_{n_2\times 1}(\mathbb{K})$,\ldots,
$[\alpha_k]\;\in\;\mathcal{M}_{n_k\times 1}(\mathbb{K})$
\begin{center}
$[U_{\sigma}]\cdot([\alpha_1]\otimes [\alpha_2]\otimes \ldots\otimes
[\alpha_k])$\;=\;$[\alpha_{\sigma(1)}]\otimes
[\alpha_{\sigma(2)}]\otimes \ldots\otimes [\alpha_{\sigma(k)}]$
\end{center}
Let $\overline{\alpha_1}$\;$\in$\;$\mathcal{E}_1$,
$\overline{\alpha_2}$\;$\in$\;$\mathcal{E}_2$, \ldots,
$\overline{\alpha_k}$\;$\in$\;$\mathcal{E}_k$ and $\mathcal{B}_1$,
$\mathcal{B}_2$,\ldots , $\mathcal{B}_k$ be basis respectively of
$\mathcal{E}_1$, $\mathcal{E}_2$, \ldots, $\mathcal{E}_k$ where the
components of $\overline{\alpha_1}$, $\overline{\alpha_2}$,\ldots,
$\overline{\alpha_k}$ form the unicolumn matrices $[\alpha_1]$,
$[\alpha_2]$,\ldots, $[\alpha_k]$.\\
$U_{\sigma}$\;$\in\;$ $\mathcal{L}$($\mathcal{E}_1$\;$\otimes$\;
$\mathcal{E}_2$\;$\otimes$\;\ldots\;$\otimes$ $\mathcal{E}_k$,
$\mathcal{E}_{\sigma(1)}$\;$\otimes$\;
$\mathcal{E}_{\sigma(2)}$\;$\otimes$\;\ldots\;$\otimes$\;
$\mathcal{E}_{\sigma(k)}$) whose matrix with respect to
($\mathcal{B}_1$$\otimes$$\mathcal{B}_2$$\otimes$ \ldots$\otimes$
$\mathcal{B}_k$, $\mathcal{B}_{\sigma(1)}$$\otimes$
$\mathcal{B}_{\sigma(2)}$$\otimes$ \ldots$\otimes$
$\mathcal{B}_{\sigma(k)}$ ) is $[U_{\sigma}]$. Thus
\begin{center}
$U_{\sigma}$($\overline{\alpha_1}$\;$\otimes$\;$\overline{\alpha_2}$\;
$\otimes$\;\ldots\;$\otimes$\;$\overline{\alpha_k}$)\;=\;
$\overline{\alpha_{\sigma(1)}}$\;$\otimes$\;$\overline{\alpha_{\sigma(2)}}$\;$\otimes$\;\ldots
\;$\otimes$\;$\overline{\alpha_{\sigma(k)}}$
\end{center}
This is true for all $\overline{\alpha_1}$\;$\in$\;$\mathcal{E}_1$,
$\overline{\alpha_2}$\;$\in$\;$\mathcal{E}_2$, \ldots,
$\overline{\alpha_k}$\;$\in$\;$\mathcal{E}_k$.\\
$U_{\sigma}$ is a $\sigma$- tensor permutation operator and
$[U_{\sigma}]$ is $\sigma$-tensor permutation matrix. \qed
\begin{thm}\label{tm42}
Let $[U_{\sigma}]$ be $\sigma$-tensor permutation matrix $n_1\otimes
n_2\otimes\ldots \otimes n_k$. Then, for all square matrices
$[A_1]$, $[A_2]$,\ldots, $[A_k]$, of dimensions respectively $n_1$,
$n_2$, \ldots, $n_k$
\begin{center}
$[U_{\sigma}]\cdot([A_1]\otimes [A_2]\otimes
\ldots\otimes[A_k]$)\;=\;$([A_{\sigma(1)}]\otimes
[A_{\sigma(2)}]\otimes \ldots\otimes
[A_{\sigma(k)}])\cdot[U_{\sigma}]$
\end{center}
\end{thm}
\pr Let $A_1$\;$\in$\; $\mathcal{L}$($\mathcal{E}_1$),
$A_2$\;$\in$\; $\mathcal{L}$($\mathcal{E}_2$), \ldots,
$A_k$\;$\in$\;$\mathcal{L}$($\mathcal{E}_k$) whose matrices with
respect to $\mathcal{B}_1$, $\mathcal{B}_2$,\ldots , $\mathcal{B}_k$
are respectively $[A_1]$, $[A_2]$,\ldots, $[A_k]$. Then
$[A_1]\otimes [A_2]\otimes \ldots\otimes[A_k]$ is the matrix of
$A_1\otimes A_2\otimes\ldots\;\otimes A_k$ with respect to
$\mathcal{B}_1$$\otimes$$\mathcal{B}_2$$\otimes$
\ldots$\otimes$ $\mathcal{B}_k$.\\
But,  $A_1\otimes A_2\otimes \ldots\otimes A_k\in\;
\mathcal{L}(\mathcal{E}_1\;\otimes\;
\mathcal{E}_2\;\otimes\;\ldots\;\otimes \mathcal{E}_k)$ and
$A_{\sigma(1)}\otimes A_{\sigma(2)}\otimes \ldots\otimes
A_{\sigma(k)}\in\;\mathcal{L}(\mathcal{E}_{\sigma(1)}\;\otimes\;
\mathcal{E}_{\sigma(2)}\;\otimes\;\ldots\;\otimes\;
\mathcal{E}_{\sigma(k)})$, thus\\
$U_{\sigma}\cdot(A_1\otimes A_2\otimes \ldots\otimes A_k$),
$(A_{\sigma(1)}\otimes A_{\sigma(2)}\otimes \ldots\otimes
A_{\sigma(k)})\cdot U_{\sigma}$$\;\in\;$
$\mathcal{L}$($\mathcal{E}_1$\;$\otimes$\;
$\mathcal{E}_2$\;$\otimes$\;\ldots\;$\otimes$ $\mathcal{E}_k$,
$\mathcal{E}_{\sigma(1)}$\;$\otimes$\;
$\mathcal{E}_{\sigma(2)}$\;$\otimes$\;\ldots\;$\otimes$\;
$\mathcal{E}_{\sigma(k)}$).\\
$[A_{\sigma(1)}]\otimes [A_{\sigma(2)}]\otimes \ldots\otimes
[A_{\sigma(k)}]$ is the matrix of $A_{\sigma(1)}\otimes
A_{\sigma(2)}\otimes \ldots\otimes A_{\sigma(k)}$ with respect to
$\mathcal{B}_{\sigma(1)}$$\otimes$
$\mathcal{B}_{\sigma(2)}$$\otimes$ \ldots$\otimes$
$\mathcal{B}_{\sigma(k)}$, thus $([A_{\sigma(1)}]\otimes
[A_{\sigma(2)}]\otimes \ldots\otimes
[A_{\sigma(k)}])\cdot[U_{\sigma}]$ is the one of
$(A_{\sigma(1)}\otimes A_{\sigma(2)}\otimes \ldots\otimes
A_{\sigma(k)})\cdot U_{\sigma}$ with respect to
($\mathcal{B}_1$$\otimes$$\mathcal{B}_2$$\otimes$ \ldots$\otimes$
$\mathcal{B}_k$, $\mathcal{B}_{\sigma(1)}$$\otimes$
$\mathcal{B}_{\sigma(2)}$$\otimes$ \ldots$\otimes$
$\mathcal{B}_{\sigma(k)}$ ).\\
 $[U_{\sigma}]\cdot([A_1]\otimes [A_2]\otimes
\ldots\otimes[A_k]$) is the matrix of $U_{\sigma}\cdot(A_1\otimes
A_2\otimes \ldots\otimes A_k$) with respect to the same basis.\\
Employing the Theorem\ref{tm32},
\begin{center}
$U_{\sigma}\cdot(A_1\otimes A_2\otimes \ldots\otimes
A_k$)\;=\;$(A_{\sigma(1)}\otimes A_{\sigma(2)}\otimes \ldots\otimes
A_{\sigma(k)})\cdot U_{\sigma}$

\end{center}
Hence
\begin{center}
$[U_{\sigma}]\cdot([A_1]\otimes [A_2]\otimes
\ldots\otimes[A_k]$)\;=\;$([A_{\sigma(1)}]\otimes
[A_{\sigma(2)}]\otimes \ldots\otimes
[A_{\sigma(k)}])\cdot[U_{\sigma}]$
\end{center}

\qed
\begin{ex}\label{tm43}
Let $[\alpha]$\;=\;$\left[%
\begin{array}{c}
  \alpha^1 \\
  \alpha^2 \\
  \alpha^3\\
\end{array}%
\right]$\;$\in$\;$\mathcal{M}_{3\times1}$($\mathbb{K}$),
$[\beta]$\;=\;$\left[%
\begin{array}{c}
  \beta^1 \\
  \beta^2 \\
\end{array}%
\right]$\;$\in$\;$\mathcal{M}_{2\times1}$($\mathbb{K}$)
\end{ex}
$[U_{3\otimes2}]$\;$\in$\;$\mathcal{M}_{6\times6}$($\mathbb{K}$) is
a tensor commutation matrix $3\otimes2$. Using the
Theorem\ref{tm41},
\begin{center}
$[U_{3\otimes2}]$$\cdot$(\;$[\alpha]$$\otimes$$[\beta]$\;)\;=\;$[\beta]$$\otimes$$[\alpha]$\\
$[U_{3\otimes2}]\cdot\left[
                        \begin{array}{c}
                          \alpha^1\beta^1 \\
                          \alpha^1\beta^2 \\
                          \alpha^2\beta^1 \\
                          \alpha^2\beta^2 \\
                          \alpha^3\beta^1 \\
                          \alpha^3\beta^2 \\
                        \end{array}
                      \right]$\;=\;$\left[
                                     \begin{array}{c}
                                       \alpha^1\beta^1 \\
                                       \alpha^2\beta^1 \\
                                       \alpha^3\beta^1 \\
                                       \alpha^1\beta^2 \\
                                       \alpha^2\beta^2 \\
                                       \alpha^3\beta^2 \\
                                     \end{array}
                                   \right]$
                      \;=\;$\left[
                                     \begin{array}{cccccc}
                                       1 & 0 & 0 & 0 & 0 & 0 \\
                                       0 & 0 & 1 & 0 & 0 & 0 \\
                                       0 & 0 & 0 & 0 & 1 & 0 \\
                                       0 & 1 & 0 & 0 & 0 & 0 \\
                                       0 & 0 & 0 & 1 & 0 & 0 \\
                                       0 & 0 & 0 & 0 & 0 & 1 \\
                                     \end{array}
                                   \right]\cdot\left[
                                                  \begin{array}{c}
                                                    \alpha^1\beta^1 \\
                                                    \alpha^1\beta^2 \\
                                                    \alpha^2\beta^1 \\
                                                    \alpha^2\beta^2 \\
                                                    \alpha^3\beta^1 \\
                                                    \alpha^3\beta^2 \\
                                                  \end{array}
                                                \right]$
\end{center}
Using Theorem\ref{tm41},
\begin{center}
 \begin{equation}\label{e41}
 [U_{3\otimes2}]\;=\;\left[\begin{array}{cccccc}
                                       1 & 0 & 0 & 0 & 0 & 0 \\
                                       0 & 0 & 1 & 0 & 0 & 0 \\
                                       0 & 0 & 0 & 0 & 1 & 0 \\
                                       0 & 1 & 0 & 0 & 0 & 0 \\
                                       0 & 0 & 0 & 1 & 0 & 0 \\
                                       0 & 0 & 0 & 0 & 0 & 1 \\
                                     \end{array}
                                   \right]
\end{equation}
\end{center}
\begin{center}

 $[U_{2\otimes3}]\;=\;\left[\begin{array}{cccccc}
                                       1 & 0 & 0 & 0 & 0 & 0 \\
                                       0 & 0 & 0 & 1 & 0 & 0 \\
                                       0 & 1 & 0 & 0 & 0 & 0 \\
                                       0 & 0 & 0 & 0 & 1 & 0 \\
                                       0 & 0 & 1 & 0 & 0 & 0 \\
                                       0 & 0 & 0 & 0 & 0 & 1 \\
                                     \end{array}
                                   \right]$

\end{center}
is the tensor commutation matrix $2\otimes3$.\\
Using the Theorem\ref{tm42}, for all
$[A]$\;$\in$\;$\mathcal{M}_{3\times3}$($\mathbb{K}$), $[B]$
$\in$\;$\mathcal{M}_{2\times2}$($\mathbb{K}$)
\begin{center}$[U_{3\otimes2}]$$\cdot$(\;$[A]$$\otimes$$[B]$\;) = (\;$[B]$$\otimes$$[A]$\;)$\cdot$$[U_{3\otimes2}]$
\end{center}
and
\begin{center}$[U_{2\otimes3}]$$\cdot$(\;$[B]$$\otimes$$[A]$\;) = (\;$[A]$$\otimes$$[B]$\;)$\cdot$$[U_{2\otimes3}]$
\end{center}
\begin{rems}
\begin{enumerate}
\item Tensor commutation matrices $1\otimes12$, $12\otimes1$, $2\otimes6$, $6\otimes2$, $3\otimes4$ et $4\otimes3$
are 12-dimensional square matrices  and are the only $12\times12$
square matrices which are tensor commutation matrices.
\item For $n\;\in\; {\mathbb{IN}}^\star$, the tensor commutation matrix $1\otimes n$ is the
$n$-dimensional unit matrix.
\item  For $n\;\in\; {\mathbb{IN}}^\star$, the tensor commutation
matrix $n\otimes n$ is a symmetrical matrix.
\item For $n\;\in\; {\mathbb{IN}}^\star$, for $p\;\in\; {\mathbb{IN}}^\star$,
tensor commutation matrices $p\otimes n$ et $n\otimes p$ are
inverses and transposes  each other.
\item If $p$ is a prime integer number , the
only  $p$-dimensional tensor commutation matrix is the p-dimensional
unit matrix, which is the tensor commutation matrix $1\otimes p$.
\item Let E be a $n$-dimensional $\mathbb{K}$-vector space, $U\;\in\; \mathcal{L}(\mathcal{E}\otimes\mathcal{E})$ a tensor
commutation operator
\begin{center}$U = U^{-1}$
\end{center}
 That is $U$ is an  involutive endomorphism.
\item For $n\;\in\; {\mathbb{IN}}^\star$, for $p\;\in\; {\mathbb{IN}}^\star$, the set
\{ $[U_{1\otimes np}]$,  $[U_{p\otimes n}]$, $[U_{n\otimes p}]$ \}
equipped of  the habitual product of matrices is  a sous-group of
the group of the permutation matrices $np\times np$.

\item We will remark \{ $[U_{1\otimes 12}]$, $[U_{2\otimes 6}]$,
$[U_{6\otimes 2}]$, $[U_{3\otimes 4}]$, $[U_{4\otimes 3}]$ \}
equipped of  the habitual product of matrices is not a sous-group of
the group of the permutation matrices $12\times 12$ because
$[U_{2\otimes 6}].[U_{4\otimes 3}]$ is not a tensor commutation
matrix.
\end{enumerate}
\end{rems}
\begin{nt}
If $U_{\sigma}$ is a tensor permutation matrix, then $\alpha
U_{\sigma}$,($\alpha\;\in\;\mathbb{K}$) is a tensor permutation
matrix too. So denote $[U_{n_1\otimes n_2\otimes\ldots\otimes
n_k}(\sigma)]$ the $\sigma$-tensor permutation matrix $n_1\otimes
n_2\otimes\ldots\otimes n_k$ formed by 0 and 1.
\end{nt}
\begin{ex}
$\sigma$ is the permutation on $\{ 1, 2, 3 \}$, such that $\sigma(1)
= 3$, $\sigma(2) = 2$, $\sigma(3) = 1$.\\
For $[\alpha]$\;=\;$\left[%
\begin{array}{c}
  \alpha^1 \\
  \alpha^2 \\
\end{array}%
\right]$,
$[\beta]$\;=\;$\left[%
\begin{array}{c}
  \beta^1 \\
  \beta^2 \\
\end{array}%
\right]$, $[\gamma]$\;=\;$\left[%
\begin{array}{c}
  \gamma^1 \\
  \gamma^2 \\
\end{array}%
\right]$\;$\in$\;$\mathcal{M}_{2\times1}$($\mathbb{K}$), using the
Theorem\ref{tm41},
\begin{center}
$[U_{2\otimes2\otimes2}(\sigma)]\cdot([\alpha]\otimes[\beta]\otimes[\gamma])\;=\;[\gamma]\otimes[\beta]\otimes[\alpha]$
\end{center}
Using the same method as the one of the Example\ref{tm43}
\end{ex}
\begin{center}
$[U_{2\otimes2\otimes2}(\sigma)]\;=\;\left[
                                       \begin{array}{cccccccc}
                                         1 & 0 & 0 & 0 & 0 & 0 & 0 & 0 \\
                                         0 & 0 & 0 & 0 & 1 & 0 & 0 & 0 \\
                                         0 & 0 & 1 & 0 & 0 & 0 & 0 & 0 \\
                                         0 & 0 & 0 & 0 & 0 & 0 & 1 & 0 \\
                                         0 & 1 & 0 & 0 & 0 & 0 & 0 & 0 \\
                                         0 & 0 & 0 & 0 & 0 & 1 & 0 & 0 \\
                                         0 & 0 & 0 & 1 & 0 & 0 & 0 & 0 \\
                                         0 & 0 & 0 & 0 & 0 & 0 & 0 & 1 \\
                                       \end{array}
                                     \right]$
\end{center}
\section{Construction of a tensor permutation matrix}
For $r\;\in\; {\mathbb{IN}}^\star$,  define
$\left[E_{ij}^{(r)}\right]$ as the matrix $r\times r$ whose elements
are zeros except the $i$-th row and $j$-th column which is equal 1.
The tensor commutation matrix $n\otimes p$ is a linear combination
of some $[E_{ij}^{(np)}]$ with coefficient 1.
\begin{center}
   \begin{equation}\label{e51}\left[E_{ij}^{(n)}\right]\otimes \left[E_{kl}^{(p)}\right]\;=\;
                                \left[E_{[p(i-1)+k][p(j-i)+l]}^{(np)}\right]
\end{equation}
\end{center}
For example, the tensor commutation matrix $[U_{3\otimes 2}]$ of
\eqref{e41}  can be written
\begin{center}
$[U_{3\otimes 2}]\;=\;\left[E_{11}^{(6)}\right] +
\left[E_{23}^{(6)}\right] + \left[E_{35}^{(6)}\right] +
\left[E_{42}^{(6)}\right] + \left[E_{54}^{(6)}\right] +
\left[E_{66}^{(6)}\right]$
\end{center}
or
\begin{center}
$[U_{3\otimes 2}]\;=\;\left[E_{11}^{(2)}\right]\otimes
\left[E_{11}^{(3)}\right]  + \left[E_{22}^{(2)}\right]\otimes
\left[E_{12}^{(3)}\right] + \left[E_{12}^{(2)}\right]\otimes
\left[E_{21}^{(3)}\right] + \left[E_{21}^{(2)}\right]\otimes
\left[E_{23}^{(3)}\right]
                      + \left[E_{11}^{(2)}\right]\otimes \left[E_{32}^{(3)}\right] +
                      \left[E_{22}^{(2)}\right]\otimes \left[E_{33}^{(3)}\right]$
\end{center}
We are going to determine the relation between $i$, $j$, $k$ and $l$
in order that the matrix $\left[E_{ij}^{(n)}\right]\otimes
\left[E_{kl}^{(p)}\right]$ make up the linear combination giving the
tensor commutation matrix $n\otimes p$. Let $[a]$ and $[b]$ be two
column matrices with $n$,
$p$ rows respectively.\\

$[a]\;=\;\left[
            \begin{array}{c}
              a^{1} \\
              \vdots \\
              a^{n} \\
            \end{array}
          \right]$
, $[b]\;=\;\left[
            \begin{array}{c}
              b^{1} \\
              \vdots \\
              b^{p} \\
            \end{array}
          \right]$

\begin{equation*}
\begin{split}
          \left[E_{ij}^{(n)}\right]\otimes
          \left[E_{kl}^{(p)}\right]\cdot([a]\otimes[b])\;
          &=\;
          \left(\left[E_{ij}^{(n)}\right]\cdot[a]\right)\otimes\left(\left[E_{kl}^{(p)}\right]\cdot[b]\right)\;\\
          &=\;\left[
                             \begin{array}{c}
                                           0 \\
                                           \vdots \\
                                           a^{j} \\
                                           \vdots \\
                                           0 \\
                                           \end{array}
                                           \right]\otimes\left[
                             \begin{array}{c}
                                           0 \\
                                           \vdots \\
                                           b^{l} \\
                                           \vdots \\
                                           0 \\
                                           \end{array}
                                           \right]\;\\
                                           &=\;\left[
                             \begin{array}{c}
                                           0 \\
                                           \vdots \\
                                           a^{j}b^{l} \\
                                           \vdots \\
                                           0 \\
                                           \end{array}
                                           \right]
                                  \end{split}
\end{equation*}

where, in these three last matrices $a^{j}$, $b^{l}$, $a^{j}b^{l}$
are
respectively at $i$-th, $k$-th, $[p(i-1)+k]$-th rows.\\
 On the other hand
\begin{center}
$[b]\otimes[a]\;=\;\left[
                             \begin{array}{c}
                                           b^{1}a^{1} \\
                                           \vdots \\
                                           b^{l}a^{j} \\
                                           \vdots \\
                                           b^{p}a^{n} \\
                                           \end{array}
                                           \right]$
\end{center}
where $b^{l}a^{j}$ is at $[n(l-1)+j]$-th row.\\

Hence
\begin{center}
\begin{equation}\label{e52}
[p(i-1)+k]\;=\;[n(l-1)+j]
\end{equation}
\end{center}
the relation between $i$, $j$, $k$, $l$ in order that
$\left[E_{ij}^{(n)}\right]\otimes
          \left[E_{kl}^{(p)}\right]$ was among the terms of the sum giving the tensor commutation
matrix $n\otimes p$. Using the relation \eqref{e51}, for each column
$\gamma$
 \begin{equation}\label{e53}
 \gamma\;=\;p(j-1) + l
\end{equation}
the element at $p(i-1) + k = [n(l-1) + j ]$-th row is equal to 1 and
the other elements of the same column are equals to zero.\\

According to \eqref{e52} and \eqref{e53} the tensor commutation
matrix $n\otimes p$ can be constructed by the following way without
doing any calculus:
\begin{rl}
Let us start in putting 1 at first row and first column, after that
let us pass into second column in going down at the rate of $n$ rows
and put 1 at this place,  then pass into  third column in going down
at the rate of $n$ rows and put 1,and so on until there is only for
us $n-1$ rows for going down (then we have obtained as number of   1
: $p$). Then pass into the next column which is the $(p + 1)$-th
column, put 1 at the second row of this column and repeat the
process until we have only $n-2$  rows for going down (then we have
obtained as number of   1 : $2p$). After that pass into the next
column which is the $(2p + 2)$-th column, put 1 at the third row of
this column and repeat the process until we have only $n - 3$ rows
for going down (then we have obtained as number of   1 : $3p$).
Continuing in this way we will have that the element at $n\times
p$-th row and $n\times p$-th column is 1.  The other elements are 0.
\end{rl}
\begin{ex}We have the tensor commutation matrix $3\otimes5$\\
 $[U_{3\otimes5}]\;=\;\left[
                      \begin{array}{ccccccccccccccc}
                        1 & 0 & 0 & 0 & 0 & 0 & 0 & 0 & 0 & 0 & 0 & 0 & 0 & 0 & 0 \\
                        0 & 0 & 0 & 0 & 0 & 1 & 0 & 0 & 0 & 0 & 0 & 0 & 0 & 0 & 0 \\
                        0 & 0 & 0 & 0 & 0 & 0 & 0 & 0 & 0 & 0 & 1 & 0 & 0 & 0 & 0 \\
                        0 & 1 & 0 & 0 & 0 & 0 & 0 & 0 & 0 & 0 & 0 & 0 & 0 & 0 & 0 \\
                        0 & 0 & 0 & 0 & 0 & 0 & 1 & 0 & 0 & 0 & 0 & 0 & 0 & 0 & 0 \\
                        0 & 0 & 0 & 0 & 0 & 0 & 0 & 0 & 0 & 0 & 0 & 1 & 0 & 0 & 0 \\
                        0 & 0 & 1 & 0 & 0 & 0 & 0 & 0& 0 & 0 & 0 & 0 & 0 & 0 & 0 \\
                        0 & 0 & 0 & 0 & 0 & 0 & 0 & 1 & 0 & 0 & 0 & 0 & 0 & 0 & 0 \\
                        0 & 0 & 0 & 0 & 0 & 0 & 0 & 0 & 0 & 0 & 0 & 0 & 1 & 0 & 0 \\
                        0 & 0 & 0 & 1 & 0 & 0 & 0 & 0 & 0 & 0 & 0 & 0 & 0 & 0 & 0 \\
                        0 & 0& 0 & 0 & 0 & 0 & 0 & 0 & 1 & 0 & 0 & 0 & 0 & 0 & 0 \\
                        0 & 0 & 0 & 0 & 0 & 0 & 0 & 0 & 0 & 0 & 0 & 0 & 0 & 1 & 0 \\
                        0 & 0 & 0 & 0 & 1 & 0 & 0 & 0 & 0 & 0 & 0 & 0 & 0 & 0 & 0 \\
                        0 & 0 & 0 & 0 & 0 & 0 & 0 & 0 & 0 & 1 & 0 & 0 & 0 & 0 & 0 \\
                        0 & 0 & 0 & 0& 0 & 0 & 0 & 0 & 0 & 0 & 0 & 0 & 0 & 0 & 1 \\
                      \end{array}
                    \right]$

\end{ex}

\begin{prop}
We have the following formula\\

 $[U_{n_1\otimes n_2\otimes\ldots\otimes
n_k}(\sigma)]\;$\\
$=\;\displaystyle\sum_{(j_{1},\ldots,j_{k})=(1,\ldots,1)}^{(n_{1},\ldots,n_{k})}
E_{[n_{k}n_{k-1}\ldots n_{3}n_{2}(i_{1}-1)+n_{k}n_{k-1}\ldots
n_{3}(i_{2}-1)+ \ldots + n_{k}(i_{k-1}-1)+i_{k}] [n_{k}\ldots
n_{3}n_{2}(j_{1}-1)+ \ldots +
n_{k}(j_{k-1}-1)+j_{k}]}^{(n_{1}n_{2}\ldots n_{k})}$

\end{prop}
\pr For any $[a_{r}]\;=\;\left(a_{r}\right)_{1\leq j_{r}\leq
n_{r}}\;\in\;\mathcal{M}_{n_{r}\times1}$,
$r\;=\;1, 2, \ldots, k$\\
 $[U_{n_1\otimes n_2\otimes\ldots\otimes
n_k}(\sigma)]\cdot\left([a_{1}]\otimes [a_{2}]\otimes\ldots\otimes
[a_{k}]\right)\;=\;[a_{\sigma(1)}]\otimes
[a_{\sigma(2)}]\otimes\ldots\otimes [a_{\sigma(k)}]$\\
 can be written
as a sum of some matrices $E_{i_{1}j_{1}}^{n_{1}}\otimes
E_{i_{2}j_{2}}^{n_{2}}\otimes \ldots \otimes
E_{i_{k}j_{k}}^{n_{k}}$.\\ Determine the conditions on $i_{1}$,
$j_{1}$, $i_{2}$, $j_{2}$, $\ldots$, $i_{k}$, $j_{k}$ so that
$E_{i_{1}j_{1}}^{n_{1}}\otimes E_{i_{2}j_{2}}^{n_{2}}\otimes \ldots
\otimes E_{i_{k}j_{k}}^{n_{k}}$ make up the sum $[U_{n_1\otimes
n_2\otimes\ldots\otimes n_k}(\sigma)]$.
\begin{center}
 $[a_{\sigma(1)}]\otimes
[a_{\sigma(2)}]\otimes\ldots\otimes [a_{\sigma(k)}]\;=\;\left(
a_{\sigma(1)}^{j_{\sigma(1)}}a_{\sigma(2)}^{j_{\sigma(2)}}\ldots
a_{\sigma(k)}^{j_{\sigma(k)}}\right)$
\end{center}
with
\begin{center}
$a_{1}^{j_{1}}a_{2}^{j_{2}}\ldots
a_{k}^{j_{k}}\;=\;a_{\sigma(1)}^{j_{\sigma(1)}}a_{\sigma(2)}^{j_{\sigma(2)}}\ldots
a_{\sigma(k)}^{j_{\sigma(k)}}$
\end{center}
is found at $[n_{\sigma(k)}n_{\sigma(k-1)}\ldots
n_{\sigma(3)}n_{\sigma(2)}(j_{\sigma(1)}-1)+ \ldots +
n_{\sigma(k)}n_{\sigma(k-1)}\ldots n_{\sigma(3)}(j_{\sigma(2)}-1) +
n_{\sigma(k)}(j_{\sigma(k-1)}-1)+j_{\sigma(k)}]$-th row.
\begin{equation*}
\begin{split}
E_{i_{1}j_{1}}^{n_{1}}\otimes E_{i_{2}j_{2}}^{n_{2}}\otimes \ldots
\otimes E_{i_{k}j_{k}}^{n_{k}}\cdot\left([a_{1}]\otimes
[a_{2}]\otimes\ldots\otimes [a_{k}]\right)\; &=\;\left[
                                                 \begin{array}{c}
                                                   0 \\
                                                   \vdots \\
                                                   a_{1}^{j_{1}} \\
                                                   \vdots \\
                                                   0 \\
                                                 \end{array}
                                               \right]\otimes\left[
                                                 \begin{array}{c}
                                                   0 \\
                                                   \vdots \\
                                                   a_{2}^{j_{2}} \\
                                                   \vdots \\
                                                   0 \\
                                                 \end{array}
                                               \right]\otimes\ldots\otimes\left[
                                                 \begin{array}{c}
                                                   0 \\
                                                   \vdots \\
                                                   a_{k}^{j_{k}} \\
                                                   \vdots \\
                                                   0 \\
                                                 \end{array}
                                               \right]\;\\
                                               &=\;
                                               \left[
                                                 \begin{array}{c}
                                                   0 \\
                                                   \vdots \\
                                                   a_{1}^{j_{1}}a_{2}^{j_{2}}\ldots
a_{k}^{j_{k}} \\
                                                   \vdots \\
                                                   0 \\
                                                 \end{array}
                                               \right]\;
                                       \end{split}
\end{equation*}
where $a_{1}^{j_{1}}$,  $a_{2}^{j_{2}}$, \ldots, $a_{k}^{j_{k}}$ and
$a_{1}^{j_{1}}a_{2}^{j_{2}}\ldots a_{k}^{j_{k}}$ are respectively at
$i_{1}$-th, $i_{2}$-th, $\ldots$-th, $i_{k}$-th and
$[n_{k}n_{k-1}\ldots n_{3}n_{2}(i_{1}-1)+n_{k}n_{k-1}\ldots
n_{3}(i_{2}-1)+ \ldots + n_{k}(i_{k-1}-1)+i_{k}]$-th rows. Thus
$E_{i_{1}j_{1}}^{n_{1}}\otimes E_{i_{2}j_{2}}^{n_{2}}\otimes \ldots
\otimes E_{i_{k}j_{k}}^{n_{k}}$ is among the terms whose sum gives
$[U_{n_1\otimes n_2\otimes\ldots\otimes n_k}(\sigma)]$ if only if\\
\begin{multline*}
n_{k}n_{k-1}\ldots n_{3}n_{2}(i_{1}-1) +n_{k}n_{k-1}\ldots
n_{3}(i_{2}-1)+ \ldots +
n_{k}(i_{k-1}-1)+i_{k}\;\\
                 =\;n_{\sigma(k)}n_{\sigma(k-1)}\ldots
n_{\sigma(3)}n_{\sigma(2)}(j_{\sigma(1)}-1)+ \ldots +
n_{\sigma(k)}n_{\sigma(k-1)}\ldots n_{\sigma(3)}(j_{\sigma(2)}-1)
\\+ n_{\sigma(k)}(j_{\sigma(k-1)}-1)+j_{\sigma(k)}.
\end{multline*}
However,\\
$E_{i_{1}j_{1}}^{n_{1}}\otimes E_{i_{2}j_{2}}^{n_{2}}\otimes \ldots
\otimes E_{i_{k}j_{k}}^{n_{k}}\;\\=\;E_{[n_{k}n_{k-1}\ldots
n_{3}n_{2}(i_{1}-1)+n_{k}n_{k-1}\ldots n_{3}(i_{2}-1)+ \ldots +
n_{k}(i_{k-1}-1)+i_{k}] [n_{k}\ldots n_{3}n_{2}(j_{1}-1)+ \ldots +
n_{k}(j_{k-1}-1)+j_{k}]}^{(n_{1}n_{2}\ldots n_{k})}$.\\
So for the column $\gamma$\\
$\gamma\;=\;n_{k}\ldots n_{3}n_{2}(j_{1}-1)+ n_{k}\ldots
n_{3}(j_{2}-1)+\ldots + n_{k}(j_{k-1}-1)+j_{k}$\\
the element at $[n_{\sigma(k)}n_{\sigma(k-1)}\ldots
n_{\sigma(3)}n_{\sigma(2)}(j_{\sigma(1)}-1)+ \ldots +
n_{\sigma(k)}n_{\sigma(k-1)}\ldots n_{\sigma(3)}(j_{\sigma(2)}-1) +
n_{\sigma(k)}(j_{\sigma(k-1)}-1)+j_{\sigma(k)}]$-th row is equal to
1 and the other elements of the same column are zeros.
 \qed
 \section{Expression of an element of a tensor permutation matrix}
 Here $n$ and $p$ are any elements of $\mathbb{N^{\star}}$. So it is a matter of
 generalizing the expression of an element of the tensor
 commutation matrix $n\otimes n$ for any $n\;\in\;\mathbb{N^{\star}}$\cite{fuj01}. At first, study
 the above example for conjecturing the expression for the more
 general case. So we follow the way in the paper \cite{fuj01}. Then write
 $[U_{3\otimes 5}]$ by the following way:\\

$[U_{3\otimes 5}]\;=\;\left[
                        \begin{array}{ccc}
                          \left[
                             \begin{array}{ccccc}
                               1 & 0 & 0 & 0 & 0 \\
                               0 & 0 & 0 & 0 & 0 \\
                               0 & 0 & 0 & 0 & 0 \\
                             \end{array}
                           \right]
                           & \left[
                             \begin{array}{ccccc}
                               0 & 0 & 0 & 0 & 0 \\
                               1 & 0 & 0 & 0 & 0 \\
                               0 & 0 & 0 & 0 & 0 \\
                             \end{array}
                           \right] & \left[
                             \begin{array}{ccccc}
                               0 & 0 & 0 & 0 & 0 \\
                               0 & 0 & 0 & 0 & 0 \\
                               1 & 0 & 0 & 0 & 0 \\
                             \end{array}
                           \right] \\
                          \left[
                             \begin{array}{ccccc}
                               0 & 1 & 0 & 0 & 0 \\
                               0 & 0 & 0 & 0 & 0 \\
                               0 & 0 & 0 & 0 & 0 \\
                             \end{array}
                           \right] & \left[
                             \begin{array}{ccccc}
                               0 & 0 & 0 & 0 & 0 \\
                               0 & 1 & 0 & 0 & 0 \\
                               0 & 0 & 0 & 0 & 0 \\
                             \end{array}
                           \right] & \left[
                             \begin{array}{ccccc}
                               0 & 0 & 0 & 0 & 0 \\
                               0 & 0 & 0 & 0 & 0 \\
                               0 & 1 & 0 & 0 & 0 \\
                             \end{array}
                           \right] \\
                          \left[
                             \begin{array}{ccccc}
                               0 & 0 & 1 & 0 & 0 \\
                               0 & 0 & 0 & 0 & 0 \\
                               0 & 0 & 0 & 0 & 0 \\
                             \end{array}
                           \right] & \left[
                             \begin{array}{ccccc}
                               0 & 0 & 0 & 0 & 0 \\
                               0 & 0 & 1 & 0 & 0 \\
                               0 & 0 & 0 & 0 & 0 \\
                             \end{array}
                           \right] & \left[
                             \begin{array}{ccccc}
                               0 & 0 & 0 & 0 & 0 \\
                               0 & 0 & 0 & 0 & 0 \\
                               0 & 0 & 1 & 0 & 0 \\
                             \end{array}
                           \right] \\
                          \left[
                             \begin{array}{ccccc}
                               0 & 0 & 0 & 1 & 0 \\
                               0 & 0 & 0 & 0 & 0 \\
                               0 & 0 & 0 & 0 & 0 \\
                             \end{array}
                           \right] & \left[
                             \begin{array}{ccccc}
                               0 & 0 & 0 & 0 & 0 \\
                               0 & 0 & 0 & 1 & 0 \\
                               0 & 0 & 0 & 0 & 0 \\
                             \end{array}
                           \right] & \left[
                             \begin{array}{ccccc}
                               0 & 0 & 0 & 0 & 0 \\
                               0 & 0 & 0 & 0 & 0 \\
                               0 & 0 & 0 & 1 & 0 \\
                             \end{array}
                           \right] \\
                          \left[
                             \begin{array}{ccccc}
                               0 & 0 & 0 & 0 & 1 \\
                               0 & 0 & 0 & 0 & 0 \\
                               0 & 0 & 0 & 0 & 0 \\
                             \end{array}
                           \right] & \left[
                             \begin{array}{ccccc}
                               0 & 0 & 0 & 0 & 0 \\
                               0 & 0 & 0 & 0 & 1 \\
                               0 & 0 & 0 & 0 & 0 \\
                             \end{array}
                           \right] & \left[
                             \begin{array}{ccccc}
                               0 & 0 & 0 & 0 & 0 \\
                               0 & 0 & 0 & 0 & 0 \\
                               0 & 0 & 0 & 0 & 1 \\
                             \end{array}
                           \right] \\
                        \end{array}
                      \right]$\\

Consider the rectangular matrices $\left[I_{n\times
p}\right]\;=\;\left(\delta_{j}^{i}\right)_{1\leq i\leq n, 1\leq j
\leq p}$, $\left[I_{p\times
n}\right]\;=\;\left(\delta_{j}^{i}\right)_{1\leq i\leq p, 1\leq j
\leq n}$, where $\delta_{j}^{i}$ is the Kronecker symbol. The matrix
\begin{center}
$\left[I_{p\times n}\right]\otimes \left[I_{n\times
p}\right]\;=\;\left(\delta_{j_{1}j_{2}}^{i_{1}i_{2}}\right)\;
=\;\left(\delta_{j_{1}}^{i_{1}}\delta_{j_{2}}^{i_{2}}\right)$
\end{center}
where,\\
$i_{1}i_{2}\;=\;11, 12, \ldots, 1n, 21, 22, \ldots, 2n, \ldots, p1,
 p2, \ldots, pn$\\
 row indices,\\
$j_{1}j_{2}\;=\;11, 12, \ldots, 1p, 21, 22, \ldots, 2p, \ldots, n1,
 n2, \ldots, np$\\
 column indices, \\
 is a $n\times p$-dimensional square matrix, which suggest us the following proposition.

\begin{prop}
\begin{equation}\label{e61}
[U_{n\otimes
p}]\;=\;\left(U_{j_{1}j_{2}}^{i_{1}i_{2}}\right)\;=\;\left(\delta_{j_{2}j_{1}}^{i_{1}i_{2}}\right)\;
=\;\left(\delta_{j_{2}}^{i_{1}}\delta_{j_{1}}^{i_{2}}\right)
\end{equation}
where,\\
 $i_{1}i_{2}\;=\;11, 12, \ldots, 1n, 21, 22, \ldots, 2n, \ldots, p1,
 p2, \ldots, pn$\\
 row indices,\\
$j_{1}j_{2}\;=\;11, 12, \ldots, 1p, 21, 22, \ldots, 2p, \ldots, n1,
 n2, \ldots, np$\\
 column indices.
\end{prop}
\pr Let $[a]\;=\;\left(a^{j_{1}}\right)_{1\leq j_{1}\leq
n}\;\in\;\mathcal{M}_{n\times1}\left(\mathbb{K}\right)$,
$[b]\;=\;\left(b^{j_{2}}\right)_{1\leq j_{2}\leq p}\;\in\;\mathcal{M}_{p\times1}\left(\mathbb{K}\right)$\\
\begin{multline*}
\left([a]\otimes [b]\right)^{i_{1}i_{2}}\longrightarrow
\left([U_{n\otimes p}]\cdot\left([a]\otimes
[b]\right)\right)^{i_{1}i_{2}}\;=\;\delta_{j_{2}}^{i_{1}}\delta_{j_{1}}^{i_{2}}a^{j_{1}}b^{j_{2}}\;
=\;\delta_{j_{2}}^{i_{1}}b^{j_{2}}\delta_{j_{1}}^{i_{2}}a^{j_{1}}\;
=\;b^{i_{1}}a^{i_{2}}\;\\
 =\;\left([b]\otimes
[a]\right)^{i_{1}i_{2}}
\end{multline*}
  \qed

  Now, we are going to generalize the formula \eqref{e61}. Let us
  consider the matrices $\left[I_{n_{\sigma(r)}\times n_{r}}\right]\;
  =\;\left(\delta_{j_{r}}^{i_{r}}\right)_{1\leq i_{r} \leq n_{\sigma(r)}, 1\leq j_{r} \leq n_{r}}$,
  $r\;=\;1, 2, \ldots,k$.

  \begin{center}
$\left[I_{n_{\sigma(1)}\times
n_{1}}\right]\otimes\left[I_{n_{\sigma(2)}\times
n_{2}}\right]\otimes\ldots\otimes\left[I_{n_{\sigma(k)}\times
n_{k}}\right]\;=\;\left(\delta_{j_{1}j_{2}\ldots
j_{k}}^{i_{1}i_{2}\ldots i_{k}}\right)\;
=\;\left(\delta_{j_{1}}^{i_{1}}\delta_{j_{2}}^{i_{2}}\ldots
\delta_{j_{k}}^{i_{k}}\right)$,
  \end{center}
  is a $n_{1}n_{2}\ldots n_{k}\times n_{\sigma(1)}n_{\sigma(2)}\ldots n_{\sigma(k)}$-dimensional
  square matrix, which suggest us the following proposition.
\begin{prop}
\begin{equation*}
[U_{n_1\otimes n_2\otimes\ldots\otimes
n_k}(\sigma)]\;=\;\left(U_{j_{1}j_{2}\ldots j_{k}}^{i_{1}i_{2}\ldots
i_{k}}\right)\;=\;\left(\delta_{j_{\sigma(1)}}^{i_{1}}\delta_{j_{\sigma(2)}}^{i_{2}}\ldots
\delta_{j_{\sigma(k)}}^{i_{k}}\right)
\end{equation*}
\end{prop}
\pr For $\left[a_{r}\right]\;=\;\left(a_{r}^{j_{r}}\right)_{1\leq
j_{r} \leq
n_{r}}\;\in\;\mathcal{M}_{n_{r}\times1}\left(\mathbb{K}\right)$,
$r\;=\;1, 2, \ldots,k$,\\
\begin{multline*}
\left(\left[a_{1}\right]\otimes \left[a_{2}\right]\otimes\ldots
\otimes\left[a_{k}\right]\right)^{i_{1}i_{2}\ldots
i_{k}}\longrightarrow \left([U_{n_1\otimes n_2\otimes\ldots\otimes
n_k}(\sigma)]\cdot\left(\left[a_{1}\right]\otimes
\left[a_{2}\right]\otimes\ldots
\otimes\left[a_{k}\right]\right)\right)^{i_{1}i_{2}\ldots
i_{k}}\;=\;\\
=\;\delta_{j_{\sigma(1)}}^{i_{1}}\delta_{j_{\sigma(2)}}^{i_{2}}\ldots
\delta_{j_{\sigma(k)}}^{i_{k}}a_{1}^{j_{1}}a_{2}^{j_{2}}\ldots
a_{k}^{j_{k}}\;\\
=\;\delta_{j_{\sigma(1)}}^{i_{1}}a_{\sigma(1)}^{j_{\sigma(1)}}
\delta_{j_{\sigma(2)}}^{i_{2}}a_{\sigma(2)}^{j_{\sigma(2)}}\ldots
\delta_{j_{\sigma(k)}}^{i_{k}}a_{\sigma(k)}^{j_{\sigma(k)}}\;\\
=\;a_{\sigma(1)}^{i_{1}}a_{\sigma(2)}^{i_{2}}\ldots
a_{\sigma(k)}^{i_{k}}\;\\
=\;\left(\left[a_{\sigma(1)}\right]\otimes
\left[a_{\sigma(2)}\right]\otimes\ldots
\otimes\left[a_{\sigma(k)}\right]\right)^{i_{1}i_{2}\ldots i_{k}}
\end{multline*}
 \qed
 \section{Expression of $\left[U_{3\otimes 3}\right]$ by using the Gell-Mann matrices}
 The Gell-Mann matrices\cite{itzub85} are\\

$\lambda_{1}=\left[
                \begin{array}{ccc}
                  0 & 1 & 0 \\
                  1 & 0 & 0 \\
                  0 & 0 & 0 \\
                \end{array}
              \right]$,\; $\lambda_{2}=\left[
                \begin{array}{ccc}
                  0 & -i & 0 \\
                  i & 0 & 0 \\
                  0 & 0 & 0 \\
                \end{array}
              \right]$,\; $\lambda_{3}=\left[
                \begin{array}{ccc}
                  1 & 0 & 0 \\
                  0 & -1 & 0 \\
                  0 & 0 & 0 \\
                \end{array}
              \right]$,\;\\

               $\lambda_{4}=\left[
                \begin{array}{ccc}
                  0 & 0 & 1 \\
                  0 & 0 & 0 \\
                  1 & 0 & 0 \\
                \end{array}
              \right]$,
               $\lambda_{5}=\left[
                \begin{array}{ccc}
                  0 & 0 & -i \\
                  0 & 0 & 0 \\
                  i & 0 & 0 \\
                \end{array}
              \right]$,\; $\lambda_{6}=\left[
                \begin{array}{ccc}
                  0 & 0 & 0 \\
                  0 & 0 & 1 \\
                  0 & 1 & 0 \\
                \end{array}
              \right]$,\; $\lambda_{7}=\left[
                \begin{array}{ccc}
                  0 & 0 & 0 \\
                  0 & 0 & -i \\
                  0 & i & 0 \\
                \end{array}
              \right]$,\; $\lambda_{8}=\frac{1}{\sqrt{3}}\left[
                \begin{array}{ccc}
                  1 & 0 & 0 \\
                  0 & 1 & 0 \\
                  0 & 0 & -2 \\
                \end{array}
              \right]$\\

$\lambda_{0}$ denotes the 3-dimensional unit matrix. $\left(\lambda
_{i}\right)_{0\leq i\leq 8}$ forms a basis of the vector spaces of
the 3-dimensional square matrices. So according to the Theorem
\ref{thm25} the system $\left(\lambda _{i}\otimes \lambda
_{j}\right)_{0\leq\; i\; \leq\; 8,\; 0\leq\; j\; \leq\; 8}$ is a
basis of the set of 9-dimensional square matrices. Then,
$\left[U_{3\otimes 3}\right]$ can be broken down into linear
combination of this system. But being inspired by the expression of
$\left[U_{2\otimes 2}\right]$    ( Cf. formula \eqref{e1}) by using
the Pauli matrices, we have directly tried to calculate the sum
$\lambda _{i}\otimes \lambda _{i}$, $i\;=\;1, 2, \ldots, 8$.\\

$\lambda _{i}\otimes \lambda _{i}\;=\;-\;\frac{2}{3}\lambda
_{0}\otimes \lambda _{0}\;+\;2\left[
                                                            \begin{array}{ccccccccc}
                                                              1 & 0 & 0 & 0 & 0 & 0 & 0 & 0 & 0 \\
                                                              0 & 0 & 0 & 1 & 0 & 0 & 0 & 0 & 0 \\
                                                              0 & 0 & 0 & 0 & 0 & 0 & 1 & 0 & 0 \\
                                                              0 & 1 & 0 & 0 & 0 & 0 & 0 & 0 & 0 \\
                                                              0 & 0 & 0 & 0 & 1 & 0 & 0 & 0 & 0 \\
                                                              0 & 0 & 0 & 0 & 0 & 0 & 0 & 1 & 0 \\
                                                              0 & 0 & 1 & 0 & 0 & 0 & 0 & 0 & 0 \\
                                                              0 & 0 & 0 & 0 & 0 & 1 & 0 & 0 & 0 \\
                                                              0 & 0 & 0 & 0 & 0 & 0 & 0 & 0 & 1 \\
                                                            \end{array}
                                                          \right]
$\\

Hence,
\begin{center}
$\left[U_{3\otimes 3}\right]\;=\;\frac{1}{3}\lambda _{0}\otimes
\lambda _{0}\;+\;\frac{1}{2}\lambda _{i}\otimes \lambda _{i}$
\end{center}

\subsection*{CONCLUSION}

We can construct a tensor commutation matrix, with or without
calculus. We can also construct a tensor permutation matrix, but
this time, by calculus, and we have the expression of an element of
such matrix. So employing these matrices, a property of tensor
product is in addition to we have already got. As the relation
between tensor commutation matrix $2\otimes 2$ with Pauli matrices
is frequently found in quantum Information theory, we hope that its
analogous, the relation between tensor commutation matrix $3\otimes
3$ with the Gell-Mann matrices, will have also its applications in
physics.

\subsection*{Acknowledgements} The author would like to thank
Hanitriarivo Rakotoson for helpful discussion, Rakotomaniraka Hary
Niela for converting the article in word file into pdf file and Fidy
Ramamonjy for the help in Latex.

\renewcommand{\bibname}{References}


\begin{thebibliography}{9}

\bibitem{rak03} RAKOTONIRINA.C, Thèse de Doctorat de Troisième Cycle de Physique
Théorique, Université d'Antananarivo, Madagascar, (2003),
unpublished.

\bibitem{wan01} WANG.R.P, arXiv: hep-ph/0107184.

\bibitem{fuj01} FUJII.K,arXiv: quant-ph/0112090, prepared for 10th Numazu
Meeting on Integral System, Noncommutative Geometry and Quantum
theory, Numazu Shizuoka Japan, 7-9 Mai 2002.

\bibitem{fad95} FADDEV.L.D, Int.J.Mod.Phys.A, Vol.10, No 13, May,1848
(1995).

\bibitem{fra02} FRANK VERSTRAETE, Thèse de Doctorat, Katholieke Universiteit
Leuven, (2002).

\bibitem{rao86} RAOELINA ANDRIAMBOLOLONA, Algèbre linéaire et Multilinéaire.
Applications, tome 1, Collection LIRA, Madagascar, (1986).

\bibitem{raoram74} RAOELINA ANDRIAMBOLOLONA and RAMIARAMANANA.D, Ann.Univ.
Madagascar, Série Sc.Nature et Math, N°11, (1974).

\bibitem{rao77} RAOELINA ANDRIAMBOLOLONA, Ann.Univ.Madagascar, Série Sc.
Nature et Math, n°14, (1977).

\bibitem{itzub85} ITZYKSON.C and J.-B.ZUBER, Quantum Field Theory, McGraw-Hill,
Singapore, (1985).
\end{thebibliography}
\end{document}